%% file: main.tex
\documentclass{ifacconf}

\usepackage{graphicx}      
\usepackage{natbib}        

\usepackage{caption}
\usepackage{subcaption}

\usepackage{enumitem}
\usepackage{threeparttable}
\usepackage{multirow}
\usepackage{booktabs}
\usepackage{balance}

\usepackage[width=.75\textwidth]{caption}

\newtheorem{definition}{Definition}

\usepackage{tikz}
\usetikzlibrary{patterns,shapes,arrows,calc}

\begin{document}
\begin{frontmatter}

\title{Teaching Predictive Control Using Specification-based Summative Assessments\thanksref{eerpapproval}\thanksref{hipeds}}

\thanks[eerpapproval]{This work has been approved by the Imperial College London Education Ethics Review Process (Study number EERP2122-029).}
\thanks[hipeds]{Supported by the EPSRC Centre for Doctoral Training in High Performance Embedded and Distributed Systems (HiPEDS, Grant Reference EP/L016796/1) and MathWorks.}
\author[ICL_EEE]{Ian McInerney}
\author[ICL_EEE,ICL_Aero]{Eric C.\ Kerrigan}

\address[ICL_EEE]{Department of Electrical \& Electronic Engineering, Imperial College London, SW7~2AZ London, UK, email: \{i.mcinerney17,e.kerrigan\}@imperial.ac.uk}
\address[ICL_Aero]{Department of Aeronautics, Imperial College London, SW7~2AZ London, UK}

\begin{abstract}
Including Model Predictive Control (MPC) in the undergraduate/graduate control curriculum is becoming vitally important due to the growing adoption of MPC in many industrial areas.
 In this paper, we present an overview of the predictive control course taught by the authors at Imperial College London between 2018 and 2021.
 We discuss how the course evolved from focusing solely on the linear MPC formulation to covering nonlinear MPC and some of its extensions.
 We also present a novel specification-based summative assessment framework, written in MATLAB, that was developed to assess the knowledge and understanding of the students in the course by tasking them with designing a controller for a real-world problem.
 The MATLAB assessment framework was designed to provide the students with the freedom to design and implement any MPC controller they wanted.
 The submitted controllers were then assessed against over 30 variations of the real-world problem to gauge student understanding of design robustness and the MPC topics from the course.
\end{abstract}


\end{frontmatter}

\section{Introduction}

Model Predictive Control (MPC) is an advanced control method that has found widespread use in many industrial areas, including aerospace, chemical processes, power electronics, and autonomous vehicles.
 The industrial success of MPC in these areas is due to several factors, including
 \begin{enumerate}[label=(\roman*)]
     \item the ability to incorporate complex systems and constraints in the controller, \label{en:mpc:cons}
     \item a simple and cross-disciplinary conceptual idea of the controller, and \label{en:mpc:simple}
     \item increased tooling and support for the design and implementation of the controller. \label{en:mpc:tooling}
 \end{enumerate}

Factors~\ref{en:mpc:cons} and~\ref{en:mpc:simple} occur because MPC controllers are written as an optimization problem that seeks to minimize an objective subject to constraints that enforce the system dynamics model and the constraints from the designer.
 Switching between different applications is then conceptually as simple as changing the dynamics model, constraints and the objective in the optimization problem.
 Factor~\ref{en:mpc:tooling} then translates this conceptual simplicity into practice, with software tools such as FORCESPRO or MATLAB's Model Predictive Control Toolbox translating the high-level optimization problem into a working controller.
 
\begin{table*}[t!]
    \centering
    \caption{Schedule of topics covered in the courses}
    \label{tab:courseContent}
    \begin{tabular}{c||c|c}
         \textbf{Week} & \textbf{2018, 2019, and 2020 Courses} & \textbf{2021 Course}  \\
         \hline
         1 & Introduction to predictive control & Introduction to predictive control \\
         2 & Predicting the future & State-space modelling and differential equation solvers \\
         3 & Unconstrained LQR RHC & Numerical/automatic differentiation \& Discretization methods \\
         \cline{3-3}
         4 & Constrained LQR RHC & Constrained LQR RHC \\
         \cline{2-2}
         5 & Soft constraints \& Setpoint tracking & Soft constraints \& Max-type costs and constraints \\
         6 & Setpoint tracking \& Disturbance rejection & Robustness \& Constraint tightening \\
         7 & Disturbance rejection & Closed-loop stability and recursive feasibility \\
         8 & Rate constraints \& Move blocking & The real-time iteration scheme \\
         \cline{2-2}
         9 & Stability \& Robustness & Move blocking \\
         10 & Final project (2018, 2019)/NMPC direct collocation (2020) & External constraint handling methods
    \end{tabular}
\end{table*}
 
The growing use of MPC in industry is occurring alongside a shift in the expectations for university control courses.
 \citet{rossiterFirstCourseFeedback2019a} reported that a recent survey asking how a first course on control should be designed showed that respondents from industry ranked topics such as optimal control and optimal state feedback in the top 10 concepts to include, which was higher than topics such as PID and lead-lag controllers.
 In addition, industry respondents felt that pedagogical techniques such as assessments focusing on concepts and demonstrating control design using authentic simulation/implementation scenarios were as important as assessments that were focused on just the mathematical concepts/theory.
 
Designing a course covering MPC is not a simple task due to the need to balance the large set of topics that could be covered (e.g.\ stochastic, robust, nonlinear, economic) and the depth of coverage (e.g.\ optimization theory, stability, problem formulations).
 Recently, \citet{faulwasserTeachingMPCWhich2021} discussed where and how such a course on MPC can fit into the curriculum, and concluded that a first course in MPC could fit into the 2nd or 3rd year of a Bachelor's curriculum and cover an introduction to the linear-quadratic problem formulation, numerical optimization, stability, and highlight application areas.
 Further topics such as nonlinear, distributed or economic MPC could then be covered in a graduate-level course along with a more in-depth discussion of the underlying theory.
 
One possible course for undergraduates was described by \citet{honcTeachingPracticingModel2016}, with the course focusing on teaching the fundamentals of linear MPC and other topics such as model derivation, controller tuning and offset-free control, with a final project of applying MPC to the control of the water level in a tank.
 A slightly more advanced course for Masters-level students was described by \citet{kellerTeachingNonlinearModel2020}, and consisted of both linear and nonlinear MPC and related topics such as optimization theory, discretization methods, and stability theory.
 At the end of the course, the students were tasked with designing a nonlinear MPC controller for a diesel engine and implementing it in real-time on a laboratory test-bench.

In this paper, we describe the predictive control course taught at Imperial College London between 2018 and 2021 to Masters-level (MSc and final year MEng) students.
 The majority of students will have completed an introductory control course on state space and transfer function methods (whose content will vary depending on their undergraduate institution), and may not have taken any prior courses on optimization.
 This course has evolved over the four years, starting with only linear-quadratic MPC and some extensions in 2018 and transitioning to nonlinear MPC and its dependencies in 2021.
 The course includes laboratory exercises utilizing a laboratory-scale gantry crane, and also uses MATLAB Grader to provide checkpointing assessments for the students during the course.
 Instead of utilizing a written/exam-based summative assessment at the end of the course, we developed a specification-based controller design problem that allows for a more thorough assessment of the student's knowledge and understanding of MPC.
 
In these summative assessments, the student is given only the specification describing a real-world problem, which is moving an overhead gantry crane to a target point in a limited time while avoiding obstacles.
 We developed a MATLAB/Simulink framework (initially described by the authors in \citet{McInerneyPredictiveControlAssessment2018}) that provides a closed-loop simulation environment where the students write their own MATLAB functions to implement the target generator, state estimator and controller.
 The student controllers are automatically tested on a set of over 30 variations of the real-world problem that are generated by changing the obstacles/constraints and adding uncertainty to the simulation model of the gantry crane.

\section{Course Lectures}
\label{sec:courseStructure}

The course was taught over a 10-week period each year, with one 2-hour long lecture session per week.
 The topics covered in the courses can be seen in Table~\ref{tab:courseContent}.

The 2018 and 2019 courses focused exclusively on linear MPC, with the 2020 course adding a 1-lecture introduction to direct collocation-based nonlinear MPC in week 10.
 The 2021 course was redesigned to primarily focus on nonlinear MPC and its associated prerequisites, and contained only a 1-lecture overview of linear MPC.

\subsection{2018, 2019 \& 2020 courses}

The main focus of the course for 2018--2020 was to teach how to apply linear predictive control to a system by using the Linear Quadratic Regulator (LQR) formulation of MPC as a Receding Horizon Controller (RHC).
 Notably, this course did not contain any lectures dedicated to optimization theory or optimization solvers, with the students instructed to instead use the \textit{mpcqpsolver} function when implementing the MPC controllers in the course.

Instead of utilizing a singular textbook for the course, the students were given a reference list containing approximately 30 books and papers that covered the material taught in the course.

\subsubsection{Course topics}

The course was split into three parts, with the basics of linear MPC in the first, several extensions to linear MPC in the second, and more advanced topics in the third.
 The first part began by introducing feedback control and the idea of predicting the future trajectory of a system.
 Then the unconstrained finite-horizon LQR formulation was introduced along with the algorithms to construct the prediction matrix and receding-horizon state feedback controller matrix.
 Finally, both the uncondensed and condensed constrained LQR formulations were introduced, along with the algorithms to construct the appropriate constraint matrices and Hessians.

In the second part, the students were introduced to more advanced linear MPC concepts taken from various research papers in the MPC field, such as soft constraints \citep{scokaertFeasibilityIssuesLinear1999}, offset-free control \citep{pannocchiaOffsetfreeMPCExplained2015}, and move blocking \citep{cagienardMoveBlockingStrategies2007}.

The final part of the course contained more advanced topics that were not contained in either the checkpointing or summative assessments.
 In the 2018/2019 courses, the only advanced topic covered was stability theory for linear MPC, based on \citet{Mayne2000_StabilitySurvey}.
 In the 2020 course, an additional advanced topic introducing direct collocation-based nonlinear MPC was added based on the tutorial paper by \citet{kellyIntroductionTrajectoryOptimization2017}.

\subsubsection{Laboratory exercises}

The course contained two physical laboratory activities that the students performed during weeks 3/4 and 5/6.
 These activities used the INTECO 3D overhead gantry crane \citep{intecoCrane} connected to a computer running MATLAB/Simulink.

During the first laboratory, the students were given an unconstrained LQR controller that moved the crane from a starting point to a target point.
 The students were asked to modify the various controller parameters (i.e.\ cost matrices and horizon length) to see how they affect the closed loop system response and gain an intuition about the tuning of the controller.

In the second laboratory, the students were given a framework with a constrained linear MPC controller.
 The students were asked to modify the cost matrices, horizon length, and constraints (i.e.\ tighten/loosen them) to gauge their effect on the controller response.
 In this laboratory, the students were also introduced to the idea of computational complexity of the controller, with the framework recording the time used by the optimization solver at each sample so that the students could see the effect the controller changes had on the computations needed.

\subsection{2021 course}

The 2021 course was redesigned to focus on nonlinear MPC and utilized the textbook by \citet{rawlingsModelPredictiveControl2020} as a reference.
 The focus of this course was still on applying MPC, however the switch to nonlinear MPC also required the introduction of more background material to the course content. 

In addition to the switch to nonlinear MPC, the COVID-19 pandemic necessitated the course be redesigned to be a virtual course with one weekly 1-hour video lecture slot.
 This was done by using a flipped-classroom approach, where the students watched video lectures on each topic before attending the weekly session, and then participated in activities/short quizzes during the weekly session.
 
\subsubsection{Course topics} 

This course consisted of two parts, the first covering the necessary background material to formulate a nonlinear MPC problem and the second covering various formulations/extensions of nonlinear MPC.
 The background material covered in the first part included numerical differential equation solvers (e.g.\ Runge-Kutta methods, collocation methods), discretization methods for nonlinear state-space equations (e.g.\ single/multiple shooting, direct collocation), and computing derivatives (i.e.\ numerical and automatic differentiation).

In the second part of the course, the students were introduced to extensions of the nonlinear MPC formulation.
 These included soft constraints, constraint tightening and robustness \citep{saltikOutlookRobustModel2018}, the real-time iteration \citep{Gros2016}, move blocking \citep{chenEfficientMoveBlocking2020}, and external constraint handling \citep{nieExternalConstraintHandling2020}.
 
\subsubsection{Laboratory exercises}

Due to the COVID-19 pandemic, the course could not use the in-person laboratory equipment, so a new virtual laboratory was developed.
 The students were given a high-fidelty simulation model of the crane in feedback with a constrained linear MPC controller and were asked to perform the same experiments as the previous-year's in-person laboratories (i.e.\ modifying controller parameters to see the closed-loop response).

This virtual laboratory was implemented inside a MATLAB Live Script that could be run using MATLAB Online, allowing for students to do the laboratory without needing to have MATLAB installed on their computers.
 The use of the Live Script also allowed the instructions and formatted equations to be included in the same file as the code, allowing the students to more easily understand the lab.

\begin{table}[t]
    \centering
    \caption{Checkpointing assessments}
    \label{tab:assignments}
    \begin{threeparttable}
    \begin{tabular}{cc}
         \textbf{Topic} & \textbf{Problem} \\
         \hline
         \hline
         \textbf{Modelling} & Gantry crane model derivation\\
         \hline
         \multirow{3}{2.5cm}{\centering\textbf{Unconstrained RHC}} & Prediction matrix construction \\
         & Cost function matrix construction \\
         & Linear RHC law construction \\
         \hline
         \multirow{5}{2.5cm}{\centering\textbf{Constrained RHC}} & Stage constraint matrix \\
         & Trajectory constraint matrix \\
         & QP constraint matrix \\
         & Receding horizon controller \\
         & Soft constraint matrices \\
         \hline
         \multirow{3}{2.5cm}{\centering\textbf{Differential equations}\tnote{1}} & Runge-Kutta Methods \\
         & Implicit/Explicit Euler \\
         & Writing an ODE45 method \\
         \hline
         \multirow{3}{2.5cm}{\centering\textbf{Linear equations}\tnote{1}} & Solution uniqueness \\
         & Solving linear systems \\
         & Least squares method \\
         \hline
         \multirow{3}{2.5cm}{\centering\textbf{Quadrature methods}\tnote{1}} & Riemann sums \\
         & Trapezoidal method \\
         & Simpson's rule \\
         \hline
         \multirow{2}{2.5cm}{\centering\textbf{Numerical\newline optimization}\tnote{1}} & fmincon w/ manual differentiation \\
         & fmincon w/ auto differentiation \\
         \hline
    \end{tabular}
    \begin{tablenotes}
    \item[1] Only in 2021 course
    \end{tablenotes}
    \end{threeparttable}
\end{table}

\begin{figure*}[t!]
    \centering
    \begin{subfigure}[b]{0.25\textwidth}
        \centering
        \resizebox{1.0\textwidth}{!}{\input{images/wedge}}
        \caption{Shape 1: Wedge (2018 \& 2019)}
        \label{fig:shape:wedge}
    \end{subfigure}%
    \hspace*{4em}
    \begin{subfigure}[b]{0.25\textwidth}
        \centering
        \resizebox{1.0\textwidth}{!}{\input{images/cornercircle}}
        \caption{Shape 2: Circles on edge (2020)}
        \label{fig:shape:edge}
    \end{subfigure}%
    \hspace*{4em}
    \begin{subfigure}[b]{0.25\textwidth}
        \centering
        \resizebox{1.0\textwidth}{!}{\input{images/areacircles}}
        \caption{Shape 3: Circles in region (2021)}
        \label{fig:shape:circles}
    \end{subfigure}%
    \caption{Sample shape the crane must stay inside. (The allowed region is in white, forbidden region in hatched red, the starting point is the filled circle and the target point is the empty square).}
    \label{fig:shape}
\end{figure*}
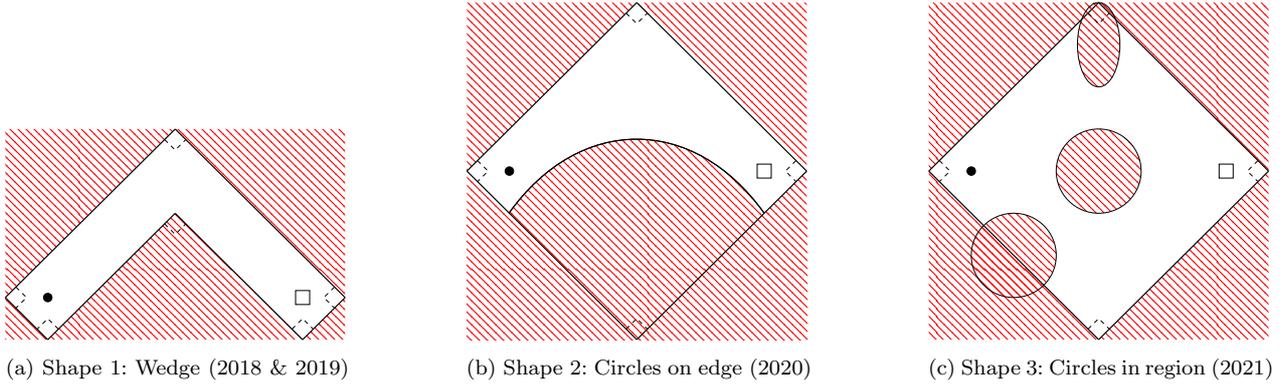

\section{Checkpointing Assessments}
\label{sec:formAssessments}

During the course, the students were assigned several small summative assessments on MATLAB Grader to gauge their progress and understanding of the topics being taught.
 In the 2018--2020 courses, there were three main topics covered in the assessments: system modelling, unconstrained RHC, and constrained RHC.
 In the 2021 course, four new topics were added to cover the new background material needed for nonlinear MPC: differential equation solvers, linear equations, quadrature methods, and numerical optimization with \textit{fmincon}.

Inside each topic, the student was presented with several MATLAB coding problems to test their knowledge on how to implement ideas discussed in the lectures (see Table~\ref{tab:assignments} for the problems in each topic).
 For example, in the prediction matrix problem in the Unconstrained RHC topic, the students were asked to write a MATLAB function that took the system's state-space matrices and desired horizon length as the input and returned the fully formed prediction matrix.
 The functions the students wrote were then tested against several sets of function inputs to ensure their function worked and was generalizable to other problems.

\section{Specification-based Summative Assessment}
\label{sec:summAssessment}

To assess the learning of the students in the course, we utilized two controller design summative assessments instead of the traditional end-of-course written exam.
 The first of the assessments was due roughly 50--60\% of the way through the course, and tested the student's knowledge of the basic predictive control concepts.
 The second assessment was due at the end of the last week of term and focused on challenging the student to try the more advanced concepts taught in part 2 of the course.
 
\subsection{Assessment Overview}

In the controller design assessments, the students were tasked with designing a controller to move a gantry crane from a starting point to a target point while staying inside a specified region and navigating around obstacles.
 For the first assessment, the constrained region was a simple rectangle (e.g.\ the left half of the wedge in Figure~\ref{fig:shape:wedge}).
 
The constrained region for the second assessment evolved during the four years of the course.
 In the 2018 and 2019 courses, the constrained region was the intersection of two rectangles at a 90$^\circ$ angle to form a wedge, as shown in Figure~\ref{fig:shape:wedge}.
 In the 2020 course, the region was changed to be a base shape of a single rectangle with circular regions centered on the edge of the rectangle added as obstacles for the crane to avoid, as shown in Figure~\ref{fig:shape:edge}.
 For the 2021 course, the shape was formed by a single rectangular region with up to 10 elliptical obstacles placed anywhere in the region, with an example region shown in Figure~\ref{fig:shape:circles}.

For the 2018--2020 courses, the objective of the controller was to simply move the crane to the target point within a specified time limit and without violating any of the constraints.
 With the switch to nonlinear MPC in the 2021 course, the assessment was updated with two major changes: (i) turning the control design problem into a fixed final time problem (i.e.\ the crane needed to be at the target point exactly at the specified time), and (ii) the addition of a constraint that limited the amount of total work the controller could perform when moving the crane.

\begin{table*}[t]
    \centering
    \caption{Number of students choosing each design option.}
    \label{tab:designchoice}
    \subfloat[2018 course (out of 30 students)\label{tab:designChoice:2018}]{
        \begin{tabular}[t]{cc}
            \toprule
            \multicolumn{2}{c}{\textbf{Cost Function}} \\
            \hline
            Quadratic cost & 30 \\
            Stabilizing terminal penalty & 15 \\
            \bottomrule
            \multicolumn{2}{c}{\textbf{Constraints}} \\
            \hline
            Soft constraints & 14 \\
            Multiple constraint sets & 29 \\
            \bottomrule
            \multicolumn{2}{c}{\textbf{State Estimator}} \\
            \hline
            Kalman filter & 15 \\
            Other state estimator & 5 \\
            \bottomrule
            \multicolumn{2}{c}{\textbf{Other Features}} \\
            \hline
            Offset-free tracking & 18 \\
            Move blocking & 1 \\
            \bottomrule
        \end{tabular}}\hspace*{2cm}
    \subfloat[2021 course (out of 26 students)\label{tab:designChoice:2021}]{
        \begin{tabular}[t]{cc}
            \toprule
            \multicolumn{2}{c}{\textbf{Setup - Path Planning}} \\
            \hline
            MATLAB \textit{nlmpc} & 2 \\
            MATLAB \textit{fmincon} & 4 \\
            Other path planning & 3 \\
            \bottomrule
            \multicolumn{2}{c}{\textbf{Setup - Nonlinearities}} \\
            \hline
            Nonlinear cost & 2 \\
            Nonlinear ellipses & 9 \\
            \bottomrule
            \multicolumn{2}{c}{\textbf{Controller - Other Features}} \\
            \hline
            Constraint tightening & 15 \\
            Soft constraints & 4 \\
            State estimator & 5 \\
            \bottomrule
        \end{tabular}
        \begin{tabular}[t]{cc}
            \toprule
            \multicolumn{2}{c}{\textbf{Controller - Optimizer}} \\
            \hline
            MATLAB \textit{nlmpc} & 4 \\
            MATLAB \textit{fmincon} & 20 \\
            Real-time iteration & 3 \\
            \bottomrule
            \multicolumn{2}{c}{\textbf{Controller - Model}} \\
            \hline
            Nonlinear/Time-varying & 9 \\
            Linear & 15 \\
            \bottomrule
            \multicolumn{2}{c}{\textbf{Controller - Cost}} \\
            \hline
            Nonlinear & 5 \\
            Quadratic & 23 \\
            \bottomrule
        \end{tabular}}
\end{table*}

\subsubsection{Specification}
 
In the assessment document, the students were only given a set of formal performance specifications that their controller must meet, and were then free to implement any MPC formulation they wished.
 These performance specifications consisted of two parts: the equilibrium condition for the crane at the target point and the definition of successful completion of a testcase, shown in Definitions~\ref{def:softEqu} and~\ref{def:softComplete}, respectively, for the 2021 course.

\begin{definition}[Equilibrium]
Obtaining equilibrium in the simulation run means that:
\begin{itemize}[noitemsep,nolistsep]
    \item the $x$ and $y$ position states of the cart are within $\epsilon_{t}$ of the target point at $t=T_f$ seconds,
    \item the $x$ and $y$ position states of the payload are within $\epsilon_{t}$ of the target point at $t=T_f$ seconds,
    \item the velocity of the cart ($\dot{x}$ and $\dot{y}$) and  angular velocity of the pendulum ($\dot{\theta}$ and $\dot{\psi})$ are within $\epsilon_{r}$ of $0$ at  $t=T_f$ seconds, and
    \item the inputs are within $\epsilon_{r}$ of $0$  at $t=T_f$ seconds.
\end{itemize}
All comparisons are made using the infinity norm.
\label{def:softEqu}
\end{definition}
\begin{definition}[Successful Completion]
Successful completion of a testcase means that:
\begin{itemize}[noitemsep,nolistsep]
    \item the system is at equilibrium (as defined in Definition~\ref{def:softEqu}) at $t=T_{f}$ seconds,
    \item the inputs remain in the interval $[-1, 1]$ during the entire time interval $[0,T_f]$,
    \item the work done by the cart over $T_f$ seconds is not more than $W_{max}$, i.e.\  $W\leq W_{max}$, and
    \item no constraints are violated over the time interval $[0,T_f]$.
\end{itemize}
\label{def:softComplete}
\end{definition}

\subsubsection{Marking}

The submitted controller designs were evaluated using both the laboratory gantry crane hardware and a high-fidelity simulation model of the gantry crane in MATLAB, and were tested against over 30 region shapes/sizes and sets of obstacles to gauge the generalizability of the student controller.

When testing against the simulation model, the students were only given one of the testcases used to mark the controllers (the default shape testcase) before the submission deadline.
 The remaining testcases were kept secret, and only given to the students along with their results.
 The secret testcases were designed to expose the controllers to a variety of situations, and were generated by
 \begin{itemize}
     \item narrowing/widening the constrained region,
     \item moving the target point closer to the constraints,
     \item adding more elliptical obstacles, and
     \item perturbing the model parameters away from the nominal model.
 \end{itemize}
 To ensure fairness in marking, every student controller was tested using the same set of secret testcases.
 
When evaluating the controller using the hardware gantry crane, each student was given a 30 minute timeslot on the actual hardware to experiment and tune their submitted controller.
 Because the hardware gantry crane already includes unmodelled elements/disturbances, no additional disturbances/uncertainties were added, and the controllers were only marked using the default shape testcase.
 
The students were required to write a short (1--2 page) report on the MPC controller they implemented.
 They also underwent an oral exam, where they were asked about their controller design and how it met the specifications given in the assessment document.

\subsection{Assessment Framework}

To administer the specification-based summative assessment, we developed a MATLAB/Simulink framework that allowed  the students to write four MATLAB functions containing their solution: \textit{mySetup}, \textit{myMPController}, \textit{myTargetGenerator}, and \textit{myStateEstimator}.
 These functions were submitted through MATLAB Grader, which also ran preliminary tests to ensure the submitted code had no syntax errors.

The marking framework  iterated through all the testcases, running a closed-loop simulation of the crane system for each testcase, and then saving the resulting state and control trajectories for later analysis.
 Each of the student functions implemented a specific component of the closed-loop system, with the \textit{mySetup} function running before the simulation began to allow offline computation of variables that were  used in the other three controller functions.

After all student controllers were run on every testcase, the marking framework then generated marks for the students by comparing the saved trajectories against the given specification and determining any violations that occurred.
 These violations were then plugged into a marking rubric to turn the testcase results into actual course marks.
 The students then received a report containing plots of all the trajectories and a listing of all the specification violations that occurred.

\section{Student Solutions and Observations}

Overall, an analysis of the student controllers submitted for the final assessment in the course shows that the specification-based summative assessment framework provided a large degree of freedom to the students.
 This can be seen in Table~\ref{tab:designchoice}, where we show the different design options chosen by students in both the 2018 and 2021 courses.

In  2018, the students were limited to using only the \textit{mpcqpsolver} optimization function (due to limitations in MATLAB's real-time code generation), meaning all 30 students implemented a quadratic MPC formulation that simply changed the constraint sets in the optimizer as the crane moved around the wedge.
 However, there was more diversity in the other design choices, with half the students choosing to implement soft constraints, another half implementing state estimation and offset-free tracking, and one student choosing to implement move blocking.
 
In the 2021 course, when there was no prescribed optimizer to use, the students chose to implement a diverse set of controllers, with the majority using \textit{fmincon}, but with 4 using the built-in MATLAB Model Predictive Control Toolbox \textit{nlmpc} function, and another 3 implementing a custom real-time iteration scheme.
 Additionally, the framework provided the freedom to implement path planning to avoid the obstacles, which 9 students chose to do (with 3 of them implementing an ASTAR-based path planner).

The majority of the students in 2021 only utilized the \textit{fmincon} optimizer to add the nonlinear elliptical constraints to the optimization problem, and still utilized a linear dynamics model and a quadratic cost function.
 To add robustness to their designs, the 2021 cohort utilized constraint tightening more than soft constraints, and relatively few implemented state estimation.
 
Based on our observations, the students engaged with and enjoyed the in-person lab components of the course, but sometimes spent too much time completing the final summative assessments.
This appears to be a result of having a lot of freedom in the design, which results in the students continually trying out new and more advanced methods (and having to spend time debugging them), only to see marginal improvements in the performance of their controller.

\section{Lessons Learned}

During the past four years, we have encountered several issues and problems when using the specification-based summative assessments in the course, with the two largest issues being: ambiguity in the specification and errors in the student code submissions.

\subsection{Students find ambiguities/loopholes}
The largest issue we faced was properly and completely defining the specification given to the students inside the assessment document so that it contained what we wanted to assess.
 When drafting the specification, we fell into the trap of including our implicit assumptions on what the controller should do when we interpreted the specification instead of interpreting it as written.
 For example, in the 2018--2020 courses the definition of equilibrium contained conditions similar to
 \begin{quote}
     The $x$ and $y$ position states of the cart are within $\epsilon_{t}$ of the target point within 5 seconds,
 \end{quote}
 which as written means the student controller only needed to be at the target point satisfying all the conditions for a single time instant to meet the specification.
 In reality, we were implicitly thinking the student controllers should enter \textit{and stay within} $\epsilon_{t}$ of the target.
 
These types of errors in the specification are the hardest to fix, since the assessment criteria could not be modified once given to the students.
 Instead, we modified the analysis framework to match the given specification, and in the 2021 course redesigned the assessment to be a fixed final time problem, removing this ambiguity.

\subsection{Don't trust student code}
While the framework was designed to be an automated system that could simply be started and then left to run, during our initial use of it we encountered many unexpected errors that crashed the framework.
 These errors mainly come from the student code, and include: incorrect computations inside the controller function causing MATLAB to throw an error, the student controller giving out-of-bound or invalid values that then cause the dynamical simulation to error, and internal errors in the numerical solvers that then crash the entire MATLAB process.
 In response, we implemented four levels of error catching/handling in the framework to catch and gracefully handle any errors that occurred, including a command line script to monitor and restart MATLAB if it crashed.



\section{Conclusions}
\balance

In this paper, we presented the predictive control course taught by the authors over the past 4 years at Imperial College London.
 This course has evolved from teaching only linear MPC in 2018 to focusing on nonlinear MPC in 2021, and assesses student knowledge of the predictive control concepts using a novel specification-based summative assessment framework.
 This framework gives the students freedom in their controller design to implement different MPC formulations, while also encouraging them to think about the robustness of their controller.
 While we observed an increase in the different types of controllers implemented by the students in the 2021 course when the focus changed to nonlinear MPC, there was still a tendency to utilize a linear model and quadratic cost instead of exploring more advanced concepts.
 Future work and updates to the course could explore ways to push students more towards the non-LQR formulations, possibly by introducing a different and more nonlinear system or changing the objective of the controller.

\bibliography{references}

\end{document}

%% file: images/wedge.tex
\begin{tikzpicture}
    \path [pattern=north west lines, pattern color=red] (1,1) -- (5,5) -- (1,5) -- (1,1);
    \path [pattern=north west lines, pattern color=red] (5,5) -- (9,1) -- (9,5) -- (5,5);
    \path [pattern=north west lines, pattern color=red] (9,1) -- (8,0) -- (9,0) -- (9,1);
    \path [pattern=north west lines, pattern color=red] (2,0) -- (5,3) -- (8,0) -- (2,0);
    \path [pattern=north west lines, pattern color=red] (1,1) -- (2,0) -- (1,0) -- (1,1);
    
    \draw [thick] (1,1) -- (5,5);
    \draw [thick] (5,5) -- (9,1);
    \draw [thick] (9,1) -- (8, 0);
    \draw [thick] (8,0) -- (5,3);
    \draw [thick] (5,3) -- (2,0);
    \draw [thick] (2,0) -- (1,1);
    
    \node [minimum size=0.5em, regular polygon, regular polygon sides=4, fill=white, draw] at (8,1) {};
    
    \draw [fill] (2,1) circle [radius=0.1];
    
    \draw [dashed] (1.25, 1.25) -- (1.5, 1.0) -- (1.25, 0.75);
    \draw [dashed] (4.75, 4.75) -- (5.0, 4.5) -- (5.25, 4.75);
    \draw [dashed] (8.75, 1.25) -- (8.5, 1.0) -- (8.75, 0.75);
    \draw [dashed] (7.75, 0.25) -- (8.0, 0.5) -- (8.25, 0.25);        
    \draw [dashed] (4.75, 2.75) -- (5.0, 2.5) -- (5.25, 2.75);        
    \draw [dashed] (2.25, 0.25) -- (2.0, 0.5) -- (1.75, 0.25);
\end{tikzpicture}

%% file: images/cornercircle.tex
\begin{tikzpicture}
    \path [pattern=north west lines, pattern color=red] (1,1) -- (5,5) -- (1,5) -- (1,1);
    \path [pattern=north west lines, pattern color=red] (5,5) -- (9,1) -- (9,5) -- (5,5);
    \path [pattern=north west lines, pattern color=red] (9,1) -- (8,0) -- (9,0) -- (9,1);
    \path [pattern=north west lines, pattern color=red] (1,1) -- (2,0) -- (1,0) -- (1,1);
    \path [pattern=north west lines, pattern color=red] (1,0) -- (1,-3) -- (9,-3) -- (9,0);
    \filldraw [pattern=north west lines, pattern color=red] (8,0) arc[start angle=45, end angle=135, x radius=4.25, y radius = 6];
    
    \draw [thick] (1,1) -- (5,5);
    \draw [thick] (5,5) -- (9,1);
    \draw [thick] (9,1) -- (8, 0);
    \draw [thick] (8,0) arc[start angle=45, end angle=135, x radius=4.25, y radius = 6];

    \draw [thick] (8, 0) -- (5, -3);
    \draw [thick] (5,-3) -- (1, 1);

    \node [minimum size=0.5em, regular polygon, regular polygon sides=4, fill=white, draw] at (8,1) {};
    
    \draw [fill] (2,1) circle [radius=0.1];
    
    \draw [dashed] (1.25, 1.25) -- (1.5, 1.0) -- (1.25, 0.75);
    \draw [dashed] (4.75, 4.75) -- (5.0, 4.5) -- (5.25, 4.75);
    \draw [dashed] (8.75, 1.25) -- (8.5, 1.0) -- (8.75, 0.75);
    \draw [dashed] (4.75, -2.75) -- (5.0, -2.5) -- (5.25, -2.75);
\end{tikzpicture}

%% file: images/areacircles.tex
\begin{tikzpicture}
    \path [pattern=north west lines, pattern color=red] (1,1) -- (5,5) -- (1,5) -- (1,1);
    \path [pattern=north west lines, pattern color=red] (5,5) -- (9,1) -- (9,5) -- (5,5);
    \path [pattern=north west lines, pattern color=red] (9,1) -- (5,-3) -- (9,-3) -- (9,1);
    \path [pattern=north west lines, pattern color=red] (1,1) -- (5,-3) -- (1,-3) -- (1,1);
    
    \draw [thick] (1,1) -- (5,5);
    \draw [thick] (5,5) -- (9,1);
    \draw [thick] (9, 1) -- (5, -3);
    \draw [thick] (5,-3) -- (1, 1);
    
    \filldraw [pattern=north west lines, pattern color=red] (5,1) circle[x radius=1, y radius = 1];

    \filldraw [pattern=north west lines, pattern color=red] (3,-1) circle[x radius=1, y radius = 1];

    \filldraw [pattern=north west lines, pattern color=red] (5,4) circle[x radius=0.5, y radius = 1];

    \node [minimum size=0.5em, regular polygon, regular polygon sides=4, fill=white, draw] at (8,1) {};
    
    \draw [fill] (2,1) circle [radius=0.1];
    
    \draw [dashed] (1.25, 1.25) -- (1.5, 1.0) -- (1.25, 0.75);
    \draw [dashed] (4.75, 4.75) -- (5.0, 4.5) -- (5.25, 4.75);
    \draw [dashed] (8.75, 1.25) -- (8.5, 1.0) -- (8.75, 0.75);
    \draw [dashed] (4.75, -2.75) -- (5.0, -2.5) -- (5.25, -2.75);
\end{tikzpicture}